\documentclass[12pt,a4paper]{article}%
\usepackage{amsmath}
\usepackage{amsfonts}
\usepackage{amssymb}
\usepackage{graphicx}
\usepackage{nopageno}%
\providecommand{\U}[1]{\protect\rule{.1in}{.1in}}
\newtheorem{theorem}{\indent{Theorem}}

\newtheorem{example}{\indent{Example}}

\newtheorem{remark}{\indent{Remark}}

\makeatletter
\renewcommand{\section}{\@startsection{section}{1}{0mm}
{\baselineskip}
{\baselineskip}{\normalfont\normalsize\scshape\centering}}
\makeatother

\begin{document}

\title{$\left.  {}\right.  $\vspace{3.5cm}\\\textbf{ON THE EXISTENCE OF EXACTLY }$N$ \textbf{LIMIT CYCLES IN LIENARD
SYSTEMS}}
\author{Aniruddha Palit$^{\ast}$ and Dhurjati Prasad Datta$^{\dagger}$\\$^{\ast}$Department of Mathematics, Surya Sen Mahavidyalaya,\\Siliguri, West Bengal, India, Pin - 734004.\\Email: mail2apalit@gmail.com\\$^{\dagger}$Department of Mathematics, University of North Bengal,\\West Bengal, India, Pin - 734013.\\Email: dp\_datta@yahoo.com}
\date{}
\maketitle

\begin{abstract}
A theorem on the existence of exactly $N$ limit cycles around a critical point
for the Lienard system $\ddot{x}+f\left(  x\right)  \dot{x}+g\left(  x\right)
=0$ is proved. An alogrithm on the determination of a desired number of limit
cycles for this system has been considered which might become relevant for a
Lienard system with incomplete data.

\bigskip

\end{abstract}

\textbf{Keywords: }Autonomous system, Lienard equation, Limit cycle

\bigskip

\textbf{AMS Classification:} 34A34,\textbf{ }70K05

\section{\textbf{Introduction}}

Limit cycles are isolated closed curves in an autonomous system in a phase
plane. Determination of shape and number of limit cycles has been a
challenging problem in the theory of autonomous systems.\ Lienard system has
been a field of active interest in recent past because of its relevance in
various physical and mathematical problem $\left[  \cite{Jordan Smith}%
-\cite{Chen Llibre Zhang}\right]  $. Recently non-smooth Lienard systems even
allowing discontinuities $\cite{Llibre Ponce Torres}$ are also being studied.

Here we consider the Lienard equation of the type%
\begin{equation}
\ddot{x}+f\left(  x\right)  \dot{x}+g\left(  x\right)  =0. \label{Lienard Eq}%
\end{equation}
The Lienard equation $\left(  \ref{Lienard Eq}\right)  $ can be written as a
non-standard autonomous system%
\begin{equation}
\dot{x}=y-F\left(  x\right)  ,\quad\dot{y}=-g\left(  x\right)
\label{Lienard System}%
\end{equation}
where $F\left(  x\right)  =\int_{0}^{x}f\left(  u\right)  du$. The phase plane
defined by $\left(  \ref{Lienard System}\right)  $ is known as Lienard plane.
Lienard gave a uniqueness theorem $\left[  \cite{Jordan Smith},\cite{Zhing
Tongren Wenzao}\right]  $ for periodic cycles for a general class of equations
when $F\left(  x\right)  $ is an odd function and satisfies a monotonicity
condition as $x\rightarrow\infty$. A challenging problem is the determination
of the number of limit cycles for a given polynomial $F\left(  x\right)  $ of
degree $\left(  m\right)  $ for the system $\left(  \ref{Lienard System}%
\right)  $ $\left[  \cite{Palit Datta},\cite{Llibre Ponce Torres},\cite{Chen
Llibre Zhang}\right]  $. Recently we have presented a new method for proving
the existence of exactly two limit cycles of a Lienard system \cite{Palit
Datta}. Recall that the proof of Lienard theorem depends on the existence of
an odd function $F\left(  x\right)  $ with zeros at $x=0$ and $x=\pm a$
$\left(  a>0\right)  $ and that $F\left(  x\right)  >0$ for $x>a$ and tends to
$\infty$ as $x\rightarrow\infty$. To weaker this assumption, we note at first
that the existence of a limit cycle is still assured if there exists a value
$\bar{\alpha}>a$ $($called an efficient upper estimate of the amplitude of the
limit cycle$)$ such that $F\left(  x\right)  $ is increasing for $a\leq
x<\bar{\alpha}<L_{1}$, where $L_{1}$ is the first extremum of $F\left(
x\right)  ,~x>a.$ Based on this observation we are then able to generalize the
standard theorem for the existence of exactly two limit cycles. Our theorem
not only extends the class of $F\left(  x\right)  $ considered by Odani
$\left[  \cite{Odani N},\cite{Odani}\right]  $, but also that of the more
recent work of Chen et al \cite{Chen Llibre Zhang} [See \cite{Palit Datta} for
more details].

In the present paper we prove the theorem for the existence of exactly $N$
limit cycles for the system $\left(  \ref{Lienard Eq}\right)  $. In the second
part of the paper we present an algorithm to generate any desired number of
limit cycles around the origin, which is the only critical point for the
system $\left(  \ref{Lienard Eq}\right)  $. Limit cycles represent an
important class of nonlinear periodic oscillations. Existence of such
nonlinear periodic cycles have been established in various natural and
biological systems $\left[  \cite{Jordan Smith},\cite{Zhing Tongren
Wenzao},\cite{Goldberger}\right]  $. It is well known that mammalian
heartbeats may follow a non-linear oscillatory patterns under certain
$($physiological$)$ constraints $\cite{Goldberger}$. However, sometimes it
becomes very difficult to obtain total information about a nonlinear system
due to various natural constraints, as a result of which we obtain only a
partial or incomplete data \cite{Donaho}. Our objective is to fill up those
gaps and construct a Lienard system that may be considered to model the
dynamics of the missing part of the phenomena in an efficient manner.

To state this in another way, let us suppose that the Lienard system is
defined only on a bounded region $\left[  -a_{1},a_{1}\right]  $, $a_{1}>0$
having one $($or at most a finite number of$)$ limit cycles in that region.
Our aim is to develop an algorithm to extend the Lienard system minimally
throughout the plane accommodating a given number of limit cycles in the
extended region. By minimal extension we mean that the graph $\left(
x,F\left(  x\right)  \right)  $, of the function $F$ which is initially
defined only in $\left\vert x\right\vert <a_{1}$ is extended beyond the line
$x=a_{1}$ iteratively as an action induced by two suitably chosen functions
$\phi\left(  x\right)  $ and $H\left(  x\right)  $ so that $\phi$ acts on the
abscissa $x$ and $H$ acts on the ordinate $F\left(  x\right)  $ respectively.
Accordingly the desired extension $\tilde{F}\left(  x\right)  $ of $F\left(
x\right)  $, $x>a_{1}$ is realized as $H\circ F\left(  x\right)  =\tilde
{F}\circ\phi\left(  x\right)  $. The choice of $\phi$ and $H$ is motivated by
theorem \ref{Th n Limit Cycle} so that the extension $\tilde{F}$ satisfies the
conditions of the said theorem. It turns out that $\phi$ can simply be a
bijective function, while $H$ may be any monotonic function admitting
$\bar{\alpha}<L$ $($c.f. equation $\left(  \ref{Definition of alfaBar}\right)
)$, $L$ being the unique extremum of $\tilde{F}\left(  x\right)  $,
$x\in\left[  a_{1},a_{2}\right]  $, $\tilde{F}\left(  a_{1}\right)  =\tilde
{F}\left(  a_{2}\right)  =0$.

The paper is organized as follows. In section \ref{Preli} we introduced our
notations. In section \ref{Existence of N limit Cycles} we have proved an
extension of the theorem in $\cite{Palit Datta}$ for existence of exactly $N$
limit cycles in the Lienard equation. In section \ref{Construction} we present
the construction by which we can get a system of the form $\left(
\ref{Lienard Eq}\right)  $ having any desired number of limit cycles around a
single critical point. Examples in support of this algorithm are studied in
section \ref{Examples}.

\section{\textbf{Notations}\label{Preli}}

We recall that \cite{Jordan Smith} by symmetry of paths, a typical phase path
$YQY^{\prime}$ of the system $\left(  \ref{Lienard System}\right)  $ becomes a
limit cycle iff $OY=OY^{\prime}$.

We consider,%
\begin{equation}
v\left(  x,y\right)  =\int_{0}^{x}g\left(  u\right)  du+\frac{1}{2}y^{2}
\label{Potential Function}%
\end{equation}
and%
\begin{equation}
v_{YQY^{\prime}}=v_{Y^{\prime}}-v_{Y}=\int\limits_{YQY^{\prime}}dv
\label{Potential Integral in v}%
\end{equation}
\begin{figure}[h]
\begin{center}
\includegraphics[height=6cm,width=4.5cm]{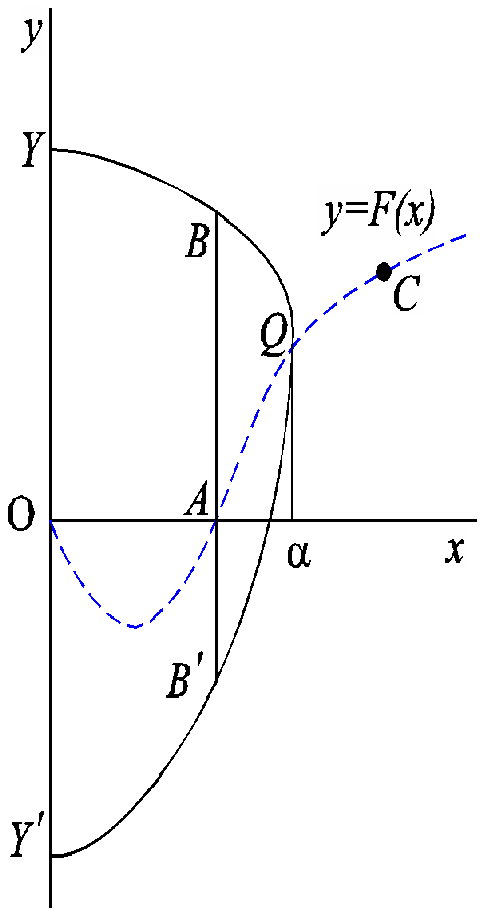}
\end{center}
\caption{{}Typical path for the Lienard theorem}%
\end{figure}It follows that%
\begin{equation}
dv=ydy+gdx=Fdy \label{Potential Differential}%
\end{equation}
so that%
\begin{equation}
OY=OY^{\prime}\Longleftrightarrow V_{YQY^{\prime}}=0\text{.}
\label{Potential Zero}%
\end{equation}
We define,%
\begin{gather*}
G\left(  x\right)  =\int_{0}^{x}g\left(  u\right)  du\\
y_{+}\left(  0\right)  =OY,\quad y_{-}\left(  0\right)  =OY^{\prime}%
\end{gather*}
and let $Q$ has coordinates $\left(  \alpha,F\left(  \alpha\right)  \right)
$. Let $\alpha^{\prime}$ and $\alpha^{\prime\prime}$ be respectively two
positive roots of the equations%
\begin{align}
G\left(  \alpha\right)   &  =\frac{1}{2}y_{+}^{2}\left(  0\right)  -\frac
{1}{2}F^{2}\left(  \alpha\right) \label{Definition of alfaDash1}\\
\text{and }G\left(  \alpha\right)   &  =\frac{1}{2}y_{-}^{2}\left(  0\right)
-\frac{1}{2}F^{2}\left(  \alpha\right)  \label{Definition of alfaDash2}%
\end{align}
Also let,%
\begin{equation}
\bar{\alpha}=\max\left\{  \alpha^{\prime},\alpha^{\prime\prime}\right\}
\label{Definition of alfaBar}%
\end{equation}
We show that $V_{YQY^{\prime}}$ has a simple zero at an $\alpha\leq\bar
{\alpha}$ \cite{Palit Datta} for the system $\left(  \ref{Lienard System}%
\right)  $ in Lienard theorem. It turns out that $\bar{\alpha}$ provides an
efficient estimate of the amplitude of the unique limit cycle of the Van der
Pol equation \cite{Palit Datta}. This result has been extended in \cite{Palit
Datta} for the existence of exactly two limit cycles as stated in the
following theorem.

\begin{theorem}
\label{Th 2 Limit Cycle}Let $f$ and $g$ be two functions satisfying the
following \linebreak properties.\newline%
\begin{tabular}
[c]{rl}%
$\left(  i\right)  $ & $f$ and $g$ are continuous;\\
$\left(  ii\right)  $ & $F$ and $g$ are odd functions and $g\left(  x\right)
>0$ for $x>0$.;\\
$\left(  iii\right)  $ & $F$ has $+ve$ simple zeros only at $x=a_{1}$,
$x=a_{2}$ for some $a_{1}>0$ and\\
& some $a_{2}>\bar{\alpha}$, $\bar{\alpha}$ being defined by $\left(
\text{\ref{Definition of alfaBar}}\right)  $ and $\bar{\alpha}<L$, where $L$
is the first\\
& local maxima of $F\left(  x\right)  $ in $\left[  a_{1},a_{2}\right]  ;$\\
$\left(  iv\right)  $ & $F$ is monotonic increasing in $a_{1}<x\leq\bar
{\alpha}$ and $F\left(  x\right)  \rightarrow-\infty$ as $x\rightarrow\infty
$\\
& monotonically for $x>a_{2};$%
\end{tabular}
\newline Then the equation $\left(  \text{\ref{Lienard Eq}}\right)  $ has
exactly two limit cycles around the origin.
\end{theorem}

It has been shown \cite{Palit Datta} that these two limit cycles are simple in
the sense that neither can bifurcate under any small $C^{1}$ perturbation
satisfying the conditions of theorem \ref{Th 2 Limit Cycle}. The existence of
$\bar{\alpha}$ satisfying an equation of the form $\left(
\ref{Definition of alfaBar}\right)  $ ensures the existence of two distinct
limit cycles.

\section{\textbf{Existence of Exactly }$N$\textbf{ limit cycles for
\protect\linebreak Lienard System}\label{Existence of N limit Cycles}}

We generalize theorem \ref{Th 2 Limit Cycle} as follows.

\begin{theorem}
\label{Th n Limit Cycle}Let $f$ and $g$ be two functions satisfying the
following \linebreak properties.\newline%
\begin{tabular}
[c]{rl}%
$\left(  i\right)  $ & $f$ and $g$ are continuous;\\
$\left(  ii\right)  $ & $F$ and $g$ are odd functions and $g\left(  x\right)
>0$ for $x>0$.;\\
$\left(  iii\right)  $ & $F$ has $N$ number of $+ve$ simple zeros only at
$x=a_{i}$, $i=1,2,\ldots,N$\\
& where $0<a_{1}<a_{2}<\ldots<a_{N}$ such that in each interval $I_{i}=\left[
a_{i},a_{i+1}\right]  $,\\
& $i=1,2,\ldots,N-1$, there exists $\bar{\alpha}_{i}$, satisfying properties
given by $\left(  \text{\ref{Definition of alfaBar}}\right)  $,\\
& such that $\bar{\alpha}_{i}<L_{i}$ where $L_{i}$ is the unique extremum in
$I_{i}$,\\
& $i=1,\ldots,N-2$ and $L_{N-1}$, the first local extremum in $\left[
a_{N-1},a_{N}\right]  $.\\
$\left(  iv\right)  $ & $F$ is monotonic in $a_{i}<x\leq\bar{\alpha}_{i}$
$\forall$ $i$ and $\left\vert F\left(  x\right)  \right\vert \rightarrow
\infty$ as $x\rightarrow\infty$\\
& monotonically for $x>a_{N}$.
\end{tabular}
\newline Then the equation $\left(  \text{\ref{Lienard Eq}}\right)  $ has
exactly $N$ limit cycles around the origin, all are simple.
\end{theorem}

\begin{figure}[h]
\begin{center}
\includegraphics[height=8cm,width=10cm]{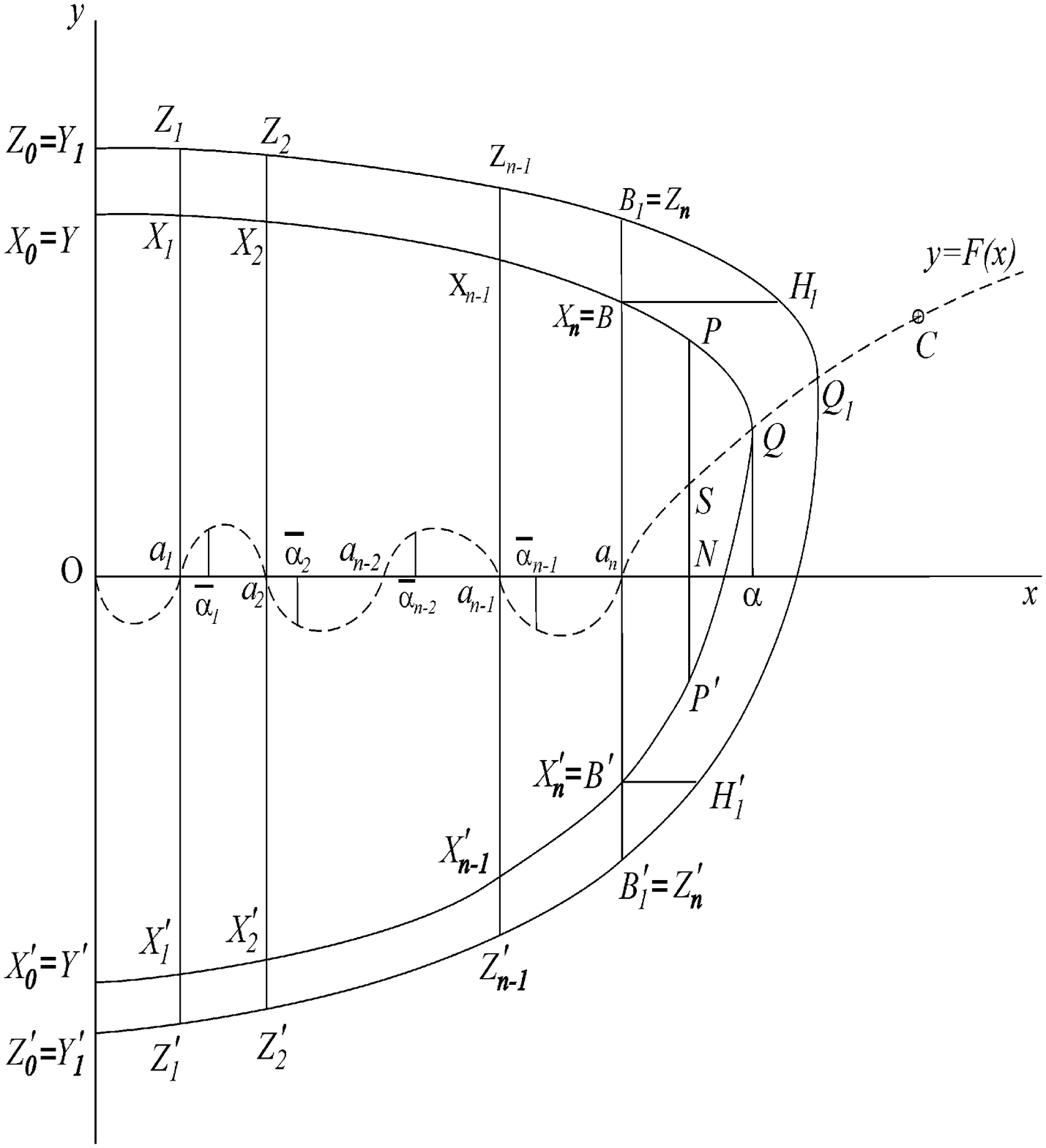}
\end{center}
\caption{{}}%
\label{Typical n Path}%
\end{figure}

\textit{Proof.} We shall prove the theorem by showing the result that each
limit cycle intersects the $x-$axis at a point lying in the open interval
$\left(  \bar{\alpha}_{i},\bar{\alpha}_{i+1}\right]  $, $i=0,1,2,\ldots,N-1$,
where $\bar{\alpha}_{0}=L_{0}$ is the local minima of $F\left(  x\right)  $ in
$\left[  0,a_{1}\right]  $. By Lienard theorem and theorem
$\ref{Th 2 Limit Cycle}$ it follows that the result is true for $N=1$ and
$N=2$. We shall now prove the theorem by induction. We assume that the theorem
is true for $N=n-1$ and we shall prove that it is true for $N=n$. We prove the
theorem by taking $n$ as an odd $+ve$ integer so that $\left(  n-1\right)  $
is even. The case for which $n$ is even can similarly be proved and so is
omitted. It can be shown that \cite{Jordan Smith}, $V_{YQY^{\prime}}$ changes
its sign from $+ve$ to $-ve$ as $Q$ moves out of $A_{1}\left(  a_{1},0\right)
$ along the curve $y=F\left(  x\right)  $ and hence vanishes there due to its
continuity and generates the first limit cycle around the origin. Next, in
\cite{Palit Datta} we see $V_{YQY^{\prime}}$ again changes its sign from $-ve$
to $+ve$ and generates the second limit cycle around the first. Also, we see
that for existence of second limit cycle we need the existence of the point
$\bar{\alpha}$, which we denote here as $\bar{\alpha}_{1}$.

Since by induction hypothesis the theorem is true for $N=n-1$, so it follows
that in each and every interval $\left(  \bar{\alpha}_{k},\bar{\alpha}%
_{k+1}\right]  $, $k=0,1,2,\ldots,n-2$ the system $\left(
\ref{Lienard System}\right)  $ has a limit cycle and the outermost limit cycle
cuts the $x-$ axis somewhere in $\left(  \bar{\alpha}_{n-1},\infty\right)  $.
Also $V_{YQY^{\prime}}$ changes its sign alternately as the point $Q$ moves
out of $a_{i}$'s, $i=1,2,\ldots,n-1$. Since $\left(  n-1\right)  $ is even, it
follows that $V_{YQY^{\prime}}$ changes its sign from $+ve$ to $-ve$ as $Q$
moves out of $a_{n-2}$ along the curve $y=F\left(  x\right)  $. Since there is
only one limit cycle in the region $\left(  \bar{\alpha}_{n-1},\infty\right)
$, so it is clear that $V_{YQY^{\prime}}$ must change its sign from $-ve$ to
$+ve$ once and only once as $Q$ moves out of $A_{n-1}\left(  a_{n-1},0\right)
$ along the curve $y=F\left(  x\right)  $. Also it follows that once
$V_{YQY^{\prime}}$ becomes $+ve$ so that it does not vanish further, otherwise
we would get one more limit cycle, contradicting the hypothesis so that total
number of limit cycle become $n$. We now try to find an estimate of $\alpha$
for which $V_{YQY^{\prime}}$ vanishes for the last time.

We shall now prove that the result is true for $N=n$ and so we assume that all
the hypotheses or conditions of this theorem are true for $N=n$. So, we get
one more point $\bar{\alpha}_{n}$ and another root $a_{n}$, ensuring the fact
that $V_{YQY^{\prime}}$ vanishes as $Q$ moves out of $A_{n-1}$ through the
curve $y=F\left(  x\right)  $, thus accommodating a unique limit cycle in the
interval $\left(  \bar{\alpha}_{n-1},\bar{\alpha}_{n}\right]  $.

By the result discussed so far it follows that $V_{YQY^{\prime}}>0$ when
$\alpha$ lies in certain suitable small right neighbourhood of $\bar{\alpha
}_{n-1}$. We shall prove that $V_{YQY^{\prime}}$ ultimately becomes $-ve$ and
remains $-ve$ as $Q$ moves out of $A_{n}\left(  a_{n},0\right)  $ along the
curve $y=F\left(  x\right)  $ generating the unique limit cycle and hence
proving the required result for $N=n$.

We draw straight line segments $X_{k}X_{k}^{\prime}$, $k=1,2,3,\ldots,n$,
passing through $A_{k}$ and parallel to $y$-axis as shown in figure
\ref{Typical n Path}. For convenience, we shall call the points $X_{n}%
,X_{n}^{\prime},Y,Y^{\prime}$ as $B,B^{\prime},X_{0},X_{0}^{\prime}$
respectively. We write the curves\newline$\left.  {}\right.  $\hfill
$\Gamma_{k}=X_{k-1}X_{k},\quad\Gamma_{k}^{\prime}=X_{k}^{\prime}%
X_{k-1}^{\prime},\quad k=1,2,3,\ldots,n$\hfill$\left.  {}\right.  $\newline so
that\newline$\left.  {}\right.  $\hfill$YQY^{\prime}=X_{0}QX_{0}^{\prime}=%
%TCIMACRO{\tsum \limits_{k=1}^{n}}%
%BeginExpansion
{\textstyle\sum\limits_{k=1}^{n}}
%EndExpansion
\Gamma_{k}+X_{n}QX_{n}^{\prime}+%
%TCIMACRO{\tsum \limits_{k=1}^{n}}%
%BeginExpansion
{\textstyle\sum\limits_{k=1}^{n}}
%EndExpansion
\Gamma_{k}^{\prime}=%
%TCIMACRO{\tsum \limits_{k=1}^{n}}%
%BeginExpansion
{\textstyle\sum\limits_{k=1}^{n}}
%EndExpansion
\left(  \Gamma_{k}+\Gamma_{k}^{\prime}\right)  +BQB^{\prime}$\hfill$\left.
{}\right.  $\newline and%
\begin{equation}
V_{YQY^{\prime}}=%
%TCIMACRO{\tsum \limits_{k=1}^{n}}%
%BeginExpansion
{\textstyle\sum\limits_{k=1}^{n}}
%EndExpansion
\left(  V_{\Gamma_{k}}+V_{\Gamma_{k}^{\prime}}\right)  +V_{BQB^{\prime}%
}\text{.} \label{Potential Sum}%
\end{equation}
We shall prove the result through the following steps.

\subsubsection*{Step $\left(  A\right)  :$ As $Q$ moves out of $A_{n}$ along
$A_{n}C$, $V_{\Gamma_{k}}+V_{\Gamma_{k}^{\prime}}$ is $+ve$ and monotonic
decreasing for odd $k$.}

$\left.  {}\right.  \hspace{18pt}$We choose two points $Q\left(
\alpha,F\left(  \alpha\right)  \right)  $ and $Q_{1}\left(  \alpha
_{1},F\left(  \alpha_{1}\right)  \right)  $ on the curve of $F\left(
x\right)  $, where $\alpha_{1}>\alpha>a_{n}$. Let $YQY^{\prime}$ and
$Y_{1}Q_{1}Y_{1}^{\prime}$ be two phase paths through $Q$ and $Q_{1}$
respectively. We have already taken $Y=X_{0}$, $Y^{\prime}=X_{0}^{\prime}$,
$B=X_{n}$ and $B^{\prime}=X_{n}^{\prime}$. We now take $Y_{1}=Z_{0}$,
$Y_{1}^{\prime}=Z_{0}^{\prime}$, $B_{1}=Z_{n}$, $B_{1}^{\prime}=Z_{n}^{\prime
}$ and $Z_{k}Z_{k}^{\prime}$ as the extension of the line segment $X_{k}%
X_{k}^{\prime}~\forall~k$. Also we write $Z_{k-1}Z_{k}=\Lambda_{k}$ and
$Z_{k}^{\prime}Z_{k-1}^{\prime}=\Lambda_{k}^{\prime}$. If $k$ is odd, then on
the segments $\Gamma_{k}$ and $\Lambda_{k}$ we have $y>0$, $F\left(  x\right)
<0$ and $y-F\left(  x\right)  >0$. Now,\newline$\left.  {}\right.  $%
\hfill$0<\left[  y-F\left(  x\right)  \right]  _{\Gamma_{k}}<\left[
y-F\left(  x\right)  \right]  _{\Lambda_{k}}$.\hfill$\left.  {}\right.
$\newline Since $g\left(  x\right)  >0$ for $x>0$ so we have\newline$\left.
{}\right.  $\hfill$\left[  \dfrac{-g\left(  x\right)  }{y-F\left(  x\right)
}\right]  _{\Gamma_{k}}<\left[  \dfrac{-g\left(  x\right)  }{y-F\left(
x\right)  }\right]  _{\Lambda_{k}}<0$.\hfill$\left.  {}\right.  $\newline So,
by $\left(  \ref{Lienard System}\right)  $ we get%
\begin{equation}
\left[  \frac{dy}{dx}\right]  _{\Gamma_{k}}<\left[  \frac{dy}{dx}\right]
_{\Lambda_{k}}<0\text{.} \label{Gradient Comparison1}%
\end{equation}
Therefore by $\left(  \ref{Gradient Comparison1}\right)  $ we have\newline%
$\left.  {}\right.  $\hfill$V_{\Gamma_{k}}=%
%TCIMACRO{\dint \limits_{\Gamma_{k}}}%
%BeginExpansion
{\displaystyle\int\limits_{\Gamma_{k}}}
%EndExpansion
F~dy=%
%TCIMACRO{\dint \limits_{\Gamma_{k}}}%
%BeginExpansion
{\displaystyle\int\limits_{\Gamma_{k}}}
%EndExpansion
\left(  -F\right)  \left(  -\dfrac{dy}{dx}\right)  dx>%
%TCIMACRO{\dint \limits_{\Lambda_{k}}}%
%BeginExpansion
{\displaystyle\int\limits_{\Lambda_{k}}}
%EndExpansion
\left(  -F\right)  \left(  -\dfrac{dy}{dx}\right)  dx=%
%TCIMACRO{\dint \limits_{\Lambda_{k}}}%
%BeginExpansion
{\displaystyle\int\limits_{\Lambda_{k}}}
%EndExpansion
F~dy=V_{\Lambda_{k}}$.\hfill$\left.  {}\right.  $\newline Since $F\left(
x\right)  $ and $dy=\dot{y}dt=-g\left(  x\right)  dt$ are both $-ve$ along
$\Lambda_{k}$ for odd $k$, so we have%
\begin{equation}
V_{\Gamma_{k}}>V_{\Lambda_{k}}=\int\limits_{\Lambda_{k}}F~dy>0\text{.}
\label{Potential Comparison1}%
\end{equation}

Next, on the segments $\Gamma_{k}^{\prime}$ and $\Lambda_{k}^{\prime}$ we have
$y<0$, $F\left(  x\right)  <0$ and $y-F\left(  x\right)  <0$. Now,\newline%
$\left.  {}\right.  $\hfill$0>\left[  y-F\left(  x\right)  \right]
_{\Gamma_{k}^{\prime}}>\left[  y-F\left(  x\right)  \right]  _{\Lambda
_{k}^{\prime}}$.\hfill$\left.  {}\right.  $\newline So, by $\left(
\ref{Lienard System}\right)  $ we get%
\begin{equation}
\left[  \frac{dy}{dx}\right]  _{\Gamma_{k}^{\prime}}>\left[  \frac{dy}%
{dx}\right]  _{\Lambda_{k}^{\prime}}>0\text{.} \label{Gradient Comparison2}%
\end{equation}
Therefore by $\left(  \ref{Gradient Comparison2}\right)  $ we have%
\[
V_{\Gamma_{k}^{\prime}}=\int\limits_{\Gamma_{k}^{\prime}}F~dy=\int
\limits_{-\Gamma_{k}^{\prime}}\left(  -F\right)  \frac{dy}{dx}dx>\int
\limits_{-\Lambda_{k}^{\prime}}\left(  -F\right)  \frac{dy}{dx}dx=\int
\limits_{\Lambda_{k}^{\prime}}F~dy=V_{\Lambda_{k}^{\prime}}\text{.}%
\]
Since $F\left(  x\right)  $ and $dy=\dot{y}dt=-g\left(  x\right)  dt$ are both
$-ve$ along $\Lambda_{k}^{\prime}$ for odd $k$, so we have%
\begin{equation}
V_{\Gamma_{k}^{\prime}}>V_{\Lambda_{k}^{\prime}}=\int\limits_{\Lambda
_{k}^{\prime}}F~dy>0\text{.} \label{Potential Comparison2}%
\end{equation}

From $\left(  \ref{Potential Comparison1}\right)  $ and $\left(
\ref{Potential Comparison2}\right)  $ we have%
\[
V_{\Gamma_{k}}+V_{\Gamma_{k}^{\prime}}>V_{\Lambda_{k}}+V_{\Lambda_{k}^{\prime
}}>0\text{.}%
\]
Therefore $V_{\Gamma_{k}}+V_{\Gamma_{k}^{\prime}}$ is $+ve$ and monotone
decreasing as the point $Q$ moves out of $A_{n}$ along $A_{n}C$.

\subsubsection*{Step $\left(  B\right)  :$ As $Q$ moves out from $A_{n}$ along
$A_{n}C$, $V_{\Gamma_{k}}+V_{\Gamma_{k}^{\prime}}$ is $-ve$ and monotonic
increasing for even $k$.}

$\left.  {}\right.  \hspace{18pt}$On the segments $\Gamma_{k}$ and
$\Lambda_{k}$ we have $y>0$, $F\left(  x\right)  >0$ and $y-F\left(  x\right)
>0$. Now,%
\[
0<\left[  y-F\left(  x\right)  \right]  _{\Gamma_{k}}<\left[  y-F\left(
x\right)  \right]  _{\Lambda_{k}}\text{.}%
\]
Since $g\left(  x\right)  >0$ for $x>0$ so we have%
\[
\left[  \frac{-g\left(  x\right)  }{y-F\left(  x\right)  }\right]
_{\Gamma_{k}}<\left[  \frac{-g\left(  x\right)  }{y-F\left(  x\right)
}\right]  _{\Lambda_{k}}<0\text{.}%
\]
So, by $\left(  \ref{Lienard System}\right)  $ we get%
\begin{equation}
\left[  \frac{dy}{dx}\right]  _{\Gamma_{k}}<\left[  \frac{dy}{dx}\right]
_{\Lambda_{k}}<0\text{.} \label{Gradient Comparison3}%
\end{equation}
Therefore by $\left(  \ref{Gradient Comparison3}\right)  $ we have%
\[
V_{\Gamma_{k}}=\int\limits_{\Gamma_{k}}F~dy=\int\limits_{\Gamma_{k}}F\frac
{dy}{dx}dx<\int\limits_{\Lambda_{k}}F\frac{dy}{dx}dx=\int\limits_{\Lambda_{k}%
}F~dy=V_{\Lambda_{k}}\text{.}%
\]
Since $F\left(  x\right)  >0$ and $dy=\dot{y}dt=-g\left(  x\right)  dt<0$
along $\Lambda_{k}$ for even $k$, so we have%
\begin{equation}
V_{\Gamma_{k}}<V_{\Lambda_{k}}=\int\limits_{\Lambda_{k}}F~dy<0\text{.}
\label{Potential Comparison3}%
\end{equation}

Next, on the segments $\Gamma_{k}^{\prime}$ and $\Lambda_{k}^{\prime}$ we have
$y<0$, $F\left(  x\right)  >0$ and $y-F\left(  x\right)  <0$. Now,%
\[
0>\left[  y-F\left(  x\right)  \right]  _{\Gamma_{k}^{\prime}}>\left[
y-F\left(  x\right)  \right]  _{\Lambda_{k}^{\prime}}%
\]
so that by $\left(  \ref{Lienard System}\right)  $ we get%
\begin{equation}
\left[  \frac{dy}{dx}\right]  _{\Gamma_{k}^{\prime}}>\left[  \frac{dy}%
{dx}\right]  _{\Lambda_{k}^{\prime}}>0\text{.} \label{Gradient Comparison4}%
\end{equation}
Therefore by $\left(  \ref{Gradient Comparison4}\right)  $ we have%
\[
V_{\Gamma_{k}^{\prime}}=\int\limits_{\Gamma_{k}^{\prime}}F~dy=\int
\limits_{-\Gamma_{k}^{\prime}}F\left(  -\frac{dy}{dx}\right)  dx<\int
\limits_{-\Lambda_{k}^{\prime}}F\left(  -\frac{dy}{dx}\right)  dx=\int
\limits_{\Lambda_{k}^{\prime}}F~dy=V_{\Lambda_{k}^{\prime}}\text{.}%
\]
Since $F\left(  x\right)  >0$ and $dy=\dot{y}dt=-g\left(  x\right)  dt<0$
along $\Lambda_{k}^{\prime}$ for even $k$, so we have%
\begin{equation}
V_{\Gamma_{k}^{\prime}}<V_{\Lambda_{k}^{\prime}}=\int\limits_{\Lambda
_{k}^{\prime}}F~dy<0\text{.} \label{Potential Comparison4}%
\end{equation}

From $\left(  \ref{Potential Comparison3}\right)  $ and $\left(
\ref{Potential Comparison4}\right)  $ we have%
\[
V_{\Gamma_{k}}+V_{\Gamma_{k}^{\prime}}<V_{\Lambda_{k}}+V_{\Lambda_{k}^{\prime
}}<0\text{.}%
\]
Therefore $V_{\Gamma_{k}}+V_{\Gamma_{k}^{\prime}}$ is $-ve$ and monotone
increasing as the point $Q$ moves out of $A_{n}$ along $A_{n}C$.

\subsubsection*{Step $\left(  C\right)  :$ $V_{BQB^{\prime}}$ is $-ve$ and
monotone decreasing and tends to $-\infty$ as $Q$ tends to infinity along
$A_{n}C$.}

$\left.  {}\right.  \hspace{18pt}$On $BQB^{\prime}$ and $B_{1}Q_{1}%
B_{1}^{\prime}$ we have $F\left(  x\right)  >0$. We draw $BH_{1}$ and
$B^{\prime}H_{1}^{\prime}$ parallel to $x$-axis. Since $F\left(  x\right)  >0$
and $dy=\dot{y}dt=-g\left(  x\right)  dt<0$ along $B_{1}Q_{1}B_{1}^{\prime}$
so%
\[
V_{B_{1}Q_{1}B_{1}^{\prime}}=\int\limits_{B_{1}Q_{1}B_{1}^{\prime}}%
F~dy\leq\int\limits_{H_{1}Q_{1}H_{1}^{\prime}}F~dy\text{.}%
\]
Since $\left[  F\left(  x\right)  \right]  _{H_{1}Q_{1}H_{1}^{\prime}}%
\geq\left[  F\left(  x\right)  \right]  _{BQB^{\prime}}$ and $F\left(
x\right)  >0$, $dy=\dot{y}dt=-g\left(  x\right)  dt<0$ along $BQB^{\prime}$ so
we have%
\begin{equation}
V_{B_{1}Q_{1}B_{1}^{\prime}}\leq\int\limits_{H_{1}Q_{1}H_{1}^{\prime}}%
F~dy\leq\int\limits_{BQB^{\prime}}F~dy=V_{BQB^{\prime}}<0\text{.}
\label{Potential Comparison5}%
\end{equation}

Let $S$ be a point on $y=F\left(  x\right)  $, to the right of $A_{n}$ and let
$BQB^{\prime}$ be an arbitrary path, with $Q$ to the right of $S$. The
straight line $PSNP^{\prime}$ is parallel to the $y$-axis. Since $F\left(
x\right)  >0$ and $dy=\dot{y}dt=-g\left(  x\right)  dt<0$ along $BQB^{\prime}$
and $PQP^{\prime}$ is a part of $BQB^{\prime}$ so we have%
\begin{equation}
V_{BQB^{\prime}}=\int\limits_{BQB^{\prime}}F~dy=-\int\limits_{B^{\prime}%
QB}F~dy\leq-\int\limits_{P^{\prime}QP}F~dy\text{.}
\label{Potential Comparison5 in y}%
\end{equation}
By hypothesis $F$ is monotone increasing on $A_{n}C$ and so $F\left(
x\right)  \geq NS$ on $PQP^{\prime}$ and hence $\left(
\ref{Potential Comparison5 in y}\right)  $ gives%
\[
V_{BQB^{\prime}}\leq-\int\limits_{P^{\prime}QP}NS\cdot dy=-NS\int
\limits_{P^{\prime}QP}dy=-NS\cdot PP^{\prime}\leq NS\cdot NP\text{.}%
\]
But as $Q$ goes to infinity towards the right, so $NP\rightarrow\infty$ and
hence by the above relation it follows that $V_{BQB^{\prime}}\rightarrow
-\infty$.

\subsubsection*{Step $\left(  D\right)  :$}

$\left.  {}\right.  \hspace{18pt}$From steps $\left(  A\right)  $ and $\left(
B\right)  $ it follows that $%
%TCIMACRO{\tsum \limits_{k=1}^{n}}%
%BeginExpansion
{\textstyle\sum\limits_{k=1}^{n}}
%EndExpansion
\left(  V_{\Gamma_{k}}+V_{\Gamma_{k}^{\prime}}\right)  $ in $\left(
\ref{Potential Sum}\right)  $ is bounded. Therefore as $Q$ moves to infinity
from the right of $A_{n}$ ultimately the quantity $V_{BQB^{\prime}}$ dominates
and hence $V_{YQY^{\prime}}$ monotonically decreases to $-\infty$ to the right
of $A_{n}$. The monotone decreasing nature of $V_{BQB^{\prime}}$ inherits the
same nature to $V_{YQY^{\prime}}$ as $Q$ moves out of $A_{n}$ along the curve
$y=F\left(  x\right)  $.

By the construction of $\bar{\alpha}_{n}$ it is clear that $V_{YQY^{\prime}%
}>0$ at a point on the left of $A_{n}$ and ultimately it becomes $-ve$ when
the point $Q$ is at the right of $A_{n}$. So, by monotonic decreasing nature
of $V_{YQY^{\prime}}$ it can vanish only once as the point $Q$ moves out of
$A_{n}$ along the curve $y=F\left(  x\right)  $. Thus, there is a unique path
for which $V_{YQY^{\prime}}=0$. By $\left(  \ref{Potential Zero}\right)  $ and
symmetry of the path it follows that the path is closed and the proof is complete.

\section{\textbf{Construction of a Lienard System with\protect\linebreak
Desired Number of Limit Cycles}\label{Construction}}

We now present an algorithm by which we can form a Lienard system with as many
limit cycles as required. We present the technique for two limit cycles around
a single critical point. This technique can similarly be extended for $n$
number of limit cycles. As stated in the introduction, \ this algorithm is
expected to become relevant in a physical model with partial or incomplete data.

Suppose in a given physical or dynamical problem, the function $F$ of the
Lienard equation $\left(  \ref{Lienard Eq}\right)  $ is well defined only with
a finite interval $\left[  -a_{1},a_{1}\right]  $ denoting $F\left(  x\right)
=f_{1}\left(  x\right)  $ for $x\in\left[  -a_{1},a_{1}\right]  $ and
satisfying the conditions:\newline%
\begin{tabular}
[c]{ll}%
$\left(  i\right)  $ & $f_{1}$ is a continuous odd function having only one
$+ve$ zero $a_{1}$\\
$\left(  ii\right)  $ & $xf_{1}\left(  x\right)  <0$~$\forall~x\in\left[
-a_{1},a_{1}\right]  $\\
$\left(  iii\right)  $ & $f_{1}$ has a unique local minimum at the point
$L_{0}$ within $\left(  a_{0},a_{1}\right)  $, $a_{0}=0$.
\end{tabular}
\newline Suppose it is also known that the system has a limit cycle just
outside the interval $\left[  -a_{1},a_{1}\right]  $. We have no information
of $F\left(  x\right)  $ beyond the interval. Our aim is to develop an
algorithm to determine a function $f_{2}$ as a restriction of $F$ in an
interval of the form $\left[  a_{1},a_{2}\right]  $ so that it satisfies the
conditions of the theorem \ref{Th n Limit Cycle}, ensuring the second limit
cycle outside $\left[  a_{1},a_{2}\right]  $. Now to determine $f_{2}$
precisely from the information of $f_{1}$ in $\left[  a_{0},a_{1}\right]  $ we
need to define two functions $\phi_{1}$ and $H_{1}$ so that we get the
abscissa and ordinates of $f_{2}$ in the interval $\left[  a_{1},a_{2}\right]
$ respectively. The choice of $\phi_{1}$ is motivated by Odani's Choice
function \cite{Odani N} $($c.f. remarks \ref{Choice Function Remark} for
further details$)$. The functions $\phi_{1}$ and $H_{1}$ are defined as
follows:\newline the functions $\phi_{1L}$ and $\phi_{1R}$ are bijective such
that%
\begin{gather*}%
\begin{array}
[c]{c}%
\phi_{1L}:\left[  a_{0},L_{0}\right]  \rightarrow\left[  a_{1},L_{1}\right]
,\quad\phi_{1L}\left(  a_{0}\right)  =a_{1},\phi_{1L}\left(  L_{0}\right)
=L_{1}\\
\phi_{1R}:\left[  L_{0},a_{1}\right]  \rightarrow\left[  L_{1},a_{2}\right]
,\quad\phi_{1R}\left(  L_{0}\right)  =L_{1},\phi_{1R}\left(  a_{1}\right)
=a_{2}%
\end{array}
\\%
\begin{array}
[c]{c}%
\phi_{1}\left(  s\right)  =\left\{
\begin{array}
[c]{c}%
\phi_{1L}\left(  s\right)  ,\quad s\in\left[  a_{0},L_{0}\right] \\
\phi_{1R}\left(  s\right)  ,\quad s\in\left[  L_{0},a_{1}\right]
\end{array}
\right.
\end{array}
\end{gather*}
and $H_{1}$ is monotone decreasing on $\left[  0,f_{1}\left(  L_{0}\right)
\right]  $ such that%
\begin{equation}
H_{1}\circ f_{1}:=f_{2}\circ\phi_{1}\text{.} \label{Objective of f2}%
\end{equation}
To make the definition $\left(  \ref{Objective of f2}\right)  $ explicit we
define at first two monotone functions $f_{2L}^{\ast}$ and $f_{2R}^{\ast}$ and
then introduce $H_{1}$ parametrically by the help of two monotone decreasing
functions $H_{1L}$ and $H_{1R}$ on $\left[  0,f_{1}\left(  L_{0}\right)
\right]  $ as%
\begin{align*}
H_{1L}  &  :f_{1L}\left(  s\right)  \rightarrow f_{2L}^{\ast}\left(  \phi
_{1L}\left(  s\right)  \right)  ,\quad s\in\left[  a_{0},L_{0}\right] \\
H_{1R}  &  :f_{1R}\left(  s\right)  \rightarrow f_{2R}^{\ast}\left(  \phi
_{1R}\left(  s\right)  \right)  ,\quad s\in\left[  L_{0},a_{1}\right] \\
H_{1}\left(  x\right)   &  =\left\{
\begin{array}
[c]{c}%
H_{1L}\left(  x\right)  ,\quad\text{if }x=f_{1L}\left(  s\right)  \text{,
}s\in\left[  a_{0},L_{0}\right] \\
H_{1R}\left(  x\right)  ,\quad\text{if }x=f_{1R}\left(  s\right)  \text{,
}s\in\left[  L_{0},a_{1}\right]
\end{array}
\right.  \text{.}%
\end{align*}
The choice of $f_{2}^{\ast}$ is made on the basis of $f_{1}\left(  x\right)  $
defined on $\left[  -a_{1},a_{1}\right]  $ and the second zero $a_{2}$ of
$F\left(  x\right)  $ that must lie close to but nevertheless, less than the
expected amplitude of the second limit cycle. We define the functions $f_{2L}$
and $f_{2R}$ as%
\begin{align*}
f_{2L}  &  :\phi_{1L}\left(  s\right)  \rightarrow H_{1L}\left(  f_{1L}\left(
s\right)  \right)  ,\quad s\in\left[  a_{0},L_{0}\right] \\
f_{2R}  &  :\phi_{1R}\left(  s\right)  \rightarrow H_{1R}\left(  f_{1R}\left(
s\right)  \right)  ,\quad s\in\left[  L_{0},a_{1}\right]  \text{.}%
\end{align*}
We should note that in the definition of $\phi_{1L}$ and $\phi_{1R}$ we have
used the conditions $\phi_{1L}\left(  a_{0}\right)  =a_{1},\phi_{1L}\left(
L_{0}\right)  =L_{1}$ and $\phi_{1R}\left(  L_{0}\right)  =L_{1},\phi
_{1R}\left(  a_{1}\right)  =a_{2}$. We could also have used the conditions
$\phi_{1L}\left(  a_{0}\right)  =L_{1},\phi_{1L}\left(  L_{0}\right)  =a_{1}$
and $\phi_{1R}\left(  L_{0}\right)  =a_{2},\phi_{1R}\left(  a_{1}\right)
=L_{1}$ instead, but in that case the function $H_{1}$ and $H_{2}$ must be
monotone increasing.\newline If $x\in\left[  a_{1},L_{1}\right]  $, then
$x=\phi_{1L}\left(  s\right)  $ for some $s\in\left[  a_{0},L_{0}\right]  $.
Therefore,%
\[
f_{2L}\left(  x\right)  =f_{2L}\left(  \phi_{1L}\left(  s\right)  \right)
=H_{1L}\left(  f_{1L}\left(  s\right)  \right)  =f_{2L}^{\ast}\left(
\phi_{1L}\left(  s\right)  \right)  =f_{2L}^{\ast}\left(  x\right)  \text{.}%
\]
So,%
\[
f_{2L}=f_{2L}^{\ast}\text{.}%
\]
Next, if $x\in\left[  L_{1},a_{2}\right]  $, then $x=\phi_{1R}\left(
s\right)  $ for some $s\in\left[  L_{0},a_{1}\right]  $. Therefore,%
\[
f_{2R}\left(  x\right)  =f_{2R}\left(  \phi_{1R}\left(  s\right)  \right)
=H_{1R}\left(  f_{1R}\left(  s\right)  \right)  =f_{2R}^{\ast}\left(
\phi_{1R}\left(  s\right)  \right)  =f_{2R}^{\ast}\left(  x\right)  \text{.}%
\]
So,%
\[
f_{2R}=f_{2R}^{\ast}\text{.}%
\]

Thus the unknown functions $f_{2L}$ and $f_{2R}$ can be expressed by known
functions $f_{2L}^{\ast}$ and $f_{2R}^{\ast}$ so that we have%
\begin{align*}
f_{2}\left(  x\right)   &  =\left\{
\begin{array}
[c]{c}%
f_{2L}\left(  x\right)  ,\quad x\in\left[  a_{1},L_{1}\right] \\
f_{2R}\left(  x\right)  ,\quad x\in\left[  L_{1},a_{2}\right]
\end{array}
\right. \\
&  =\left\{
\begin{array}
[c]{c}%
f_{2L}^{\ast}\left(  x\right)  ,\quad x\in\left[  a_{1},L_{1}\right] \\
f_{2R}^{\ast}\left(  x\right)  ,\quad x\in\left[  L_{1},a_{2}\right]
\end{array}
\right.
\end{align*}

Next we construct the restriction $f_{3}$ of the function $F$ in $\left[
a_{2},a_{3}\right]  $ having unique local maximum at $L_{2}$ $($say$)$ in
$\left(  a_{2},a_{3}\right)  $. We assume two bijective functions%
\begin{align*}
\phi_{2L}  &  :\left[  a_{1},L_{1}\right]  \rightarrow\left[  a_{2}%
,L_{2}\right]  ,\quad\phi_{2L}\left(  a_{1}\right)  =a_{2},\phi_{2L}\left(
L_{1}\right)  =L_{2}\\
\text{and }\phi_{2R}  &  :\left[  L_{1},a_{2}\right]  \rightarrow\left[
L_{2},a_{3}\right]  ,\quad\phi_{2R}\left(  L_{1}\right)  =L_{2},\phi
_{2R}\left(  a_{2}\right)  =a_{3}%
\end{align*}
and two more functions $f_{3L}^{\ast}$ and $f_{3R}^{\ast}$. We define two
monotone decreasing functions $H_{2L}$ and $H_{2R}$ on $\left[  0,f_{2}\left(
L_{1}\right)  \right]  $ parametrically as%
\begin{align*}
H_{2L}  &  :f_{2L}\left(  s\right)  \rightarrow f_{3L}^{\ast}\left(  \phi
_{2L}\left(  s\right)  \right)  ,\quad s\in\left[  a_{1},L_{1}\right] \\
H_{2R}  &  :f_{2R}\left(  s\right)  \rightarrow f_{3R}^{\ast}\left(  \phi
_{2R}\left(  s\right)  \right)  ,\quad s\in\left[  L_{1},a_{2}\right]
\text{.}%
\end{align*}
We define%
\begin{align*}
f_{3L}  &  :\phi_{2L}\left(  s\right)  \rightarrow H_{2L}\left(  f_{2L}\left(
s\right)  \right)  ,\quad s\in\left[  a_{1},L_{1}\right] \\
f_{3R}  &  :\phi_{2R}\left(  s\right)  \rightarrow H_{2R}\left(  f_{2R}\left(
s\right)  \right)  ,\quad s\in\left[  L_{1},a_{2}\right]
\end{align*}
so that as shown above we have%
\[
f_{3L}=f_{3L}^{\ast}\text{ and }f_{3R}=f_{3R}^{\ast}\text{.}%
\]
Therefore,%
\begin{align*}
f_{3}\left(  x\right)   &  =\left\{
\begin{array}
[c]{c}%
f_{3L}\left(  x\right)  ,\quad x\in\left[  a_{2},L_{2}\right] \\
f_{3R}\left(  x\right)  ,\quad x\in\left[  L_{2},a_{3}\right]
\end{array}
\right. \\
&  =\left\{
\begin{array}
[c]{c}%
f_{3L}^{\ast}\left(  x\right)  ,\quad x\in\left[  a_{2},L_{2}\right] \\
f_{3R}^{\ast}\left(  x\right)  ,\quad x\in\left[  L_{2},a_{3}\right]
\end{array}
\right.
\end{align*}
We observe that%
\begin{align*}
f_{3L}  &  :\phi_{2L}\left(  \phi_{1L}\left(  s\right)  \right)  \rightarrow
H_{2L}\left(  f_{2L}\left(  \phi_{1L}\left(  s\right)  \right)  \right)
,\quad s\in\left[  a_{0},L_{0}\right] \\
f_{3R}  &  :\phi_{2R}\left(  \phi_{1R}\left(  s\right)  \right)  \rightarrow
H_{2R}\left(  f_{2R}\left(  \phi_{1R}\left(  s\right)  \right)  \right)
,\quad s\in\left[  L_{0},a_{1}\right]  \text{.}%
\end{align*}
We can similarly proceed and construct all the restrictions $f_{k}$ of the
function $F$ in $\left[  a_{k-1},a_{k}\right]  $ for $k=4,5,6,\ldots,N$ so
that the corresponding Lienard system have exactly $N$ limit cycles. Thus an
incomplete Lienard system can be extended iteratively over larger and larger
intervals of $x$, having as many $($simple$)$ limit cycles as desired. We
note, however, that the choice of iterated functions has as large
arbitrariness except for the required minimal conditions of monotonicity
satisfying theorem \ref{Th n Limit Cycle}. The number of limit cycles for each
such choices remain invariant. The problem of reconstructing data with a given
number of limit cycles and having specified shapes is left for future study.
We now illustrate the above construction by the following examples.

\section{\textbf{Examples}\label{Examples}}

We now present some examples following the construction described in section
\ref{Construction}. Here the figures has been drawn using Mathematica 5.1.

\begin{example}
\label{Ex Construction 1}Let $a_{1}=0.2$, $a_{2}=0.5$ and%
\[
f_{1}\left(  x\right)  =0.15-0.25\sqrt{1-\frac{\left(  x-0.1\right)  ^{2}%
}{0.125^{2}}},\quad-0.2\leq x\leq0.2\text{.}%
\]
Here, $L_{0}=0.1$. Let $L_{1}=0.3$. Let us choose%
\begin{align*}
f_{2L}^{\ast}\left(  x\right)   &  =-0.15+0.25\sqrt{1-\frac{\left(
x-0.3\right)  ^{2}}{\left(  0.125\right)  ^{2}}}\\
f_{2R}^{\ast}\left(  x\right)   &  =-0.15+0.25\sqrt{1-\frac{\left(
x-0.3\right)  ^{2}}{\left(  0.25\right)  ^{2}}}\text{.}%
\end{align*}
Also, let%
\[
\phi_{1L}\left(  s\right)  =\sqrt{As^{2}+B}\text{.}%
\]
To determine the unknown parameters $A$ and $B$ we assume that $\phi
_{1L}\left(  a_{0}\right)  =a_{1}$, $\phi_{1L}\left(  L_{0}\right)  =L_{1}$.
Then $A=5$ and $B=0.04$. Next, let%
\[
\phi_{1R}\left(  s\right)  =\sqrt{A^{\prime}s^{2}+B^{\prime}}%
\]
with $\phi_{1R}\left(  L_{0}\right)  =L_{1}$ and $\phi_{1R}\left(
a_{1}\right)  =a_{2}$. Then, $A^{\prime}=\dfrac{16}{3}$ and $B^{\prime}%
=\dfrac{11}{300}$. Then following the algorithm in section $\ref{Construction}%
$ we have%
\begin{align*}
f_{2L} &  =f_{2L}^{\ast}\text{ in }\left[  a_{1},L_{1}\right]  \\
\text{and }f_{2R} &  =f_{2R}^{\ast}\text{ in }\left[  L_{1},a_{2}\right]
\end{align*}
so that%
\[
f_{2}\left(  x\right)  =\left\{
\begin{array}
[c]{c}%
f_{2L}\left(  x\right)  ,\quad x\in\left[  a_{1},L_{1}\right]  \\
f_{2R}\left(  x\right)  ,\quad x\in\left[  L_{1},a_{2}\right]
\end{array}
\right.
\]
We now define%
\[
F_{+}\left(  x\right)  =\left\{
\begin{array}
[c]{l}%
f_{1}\left(  x\right)  ,\quad0\leq x<a_{1}\\
f_{2}\left(  x\right)  ,\quad a_{1}\leq x<a_{2}\\
-\dfrac{4}{3}\left(  x-0.5\right)  ,\quad x\geq a_{2}%
\end{array}
\right.
\]
to make $F_{+}$ continuously differentiable in $\left[  0,\infty\right)  $.
The last part of the function $F_{+}$ is taken to make $F_{+}$ monotone
decreasing for $x\geq a_{2}$ so that the function $F$ defined below satisfy
the condition that $\left\vert F\left(  x\right)  \right\vert \rightarrow
\infty$ as $x\rightarrow\infty$ monotonically for $x\geq a_{2}$. We take%
\[
F\left(  x\right)  =\left\{
\begin{array}
[c]{c}%
F_{+}\left(  x\right)  ,\quad x\geq0\\
F_{-}\left(  x\right)  ,\quad x<0
\end{array}
\right.
\]
We find two limit cycles which cross the $+ve$ $y$-axis at the points $\left(
0,0.26731065\right)  $ and $\left(  0,0.5749823\right)  $ respectively. So,
$y_{+}\left(  0\right)  =y_{-}\left(  0\right)  =0.26731065$ and $\bar{\alpha
}_{1}=0.254219124$. So, the conditions $\bar{\alpha}_{1}\leq L_{1}$ are
satisfied in this example. Thus the existence of limit cycles are ensured by
the theorem $\ref{Th n Limit Cycle}$ with $g\left(  x\right)  =x$ establishing
the construction in section $\ref{Construction}$. The limit cycles alongwith
the curve of $F\left(  x\right)  $ has been shown in figure
$\ref{Ex Construction 1 Fig}$.\begin{figure}[h]
\begin{center}
\includegraphics[height=6cm]{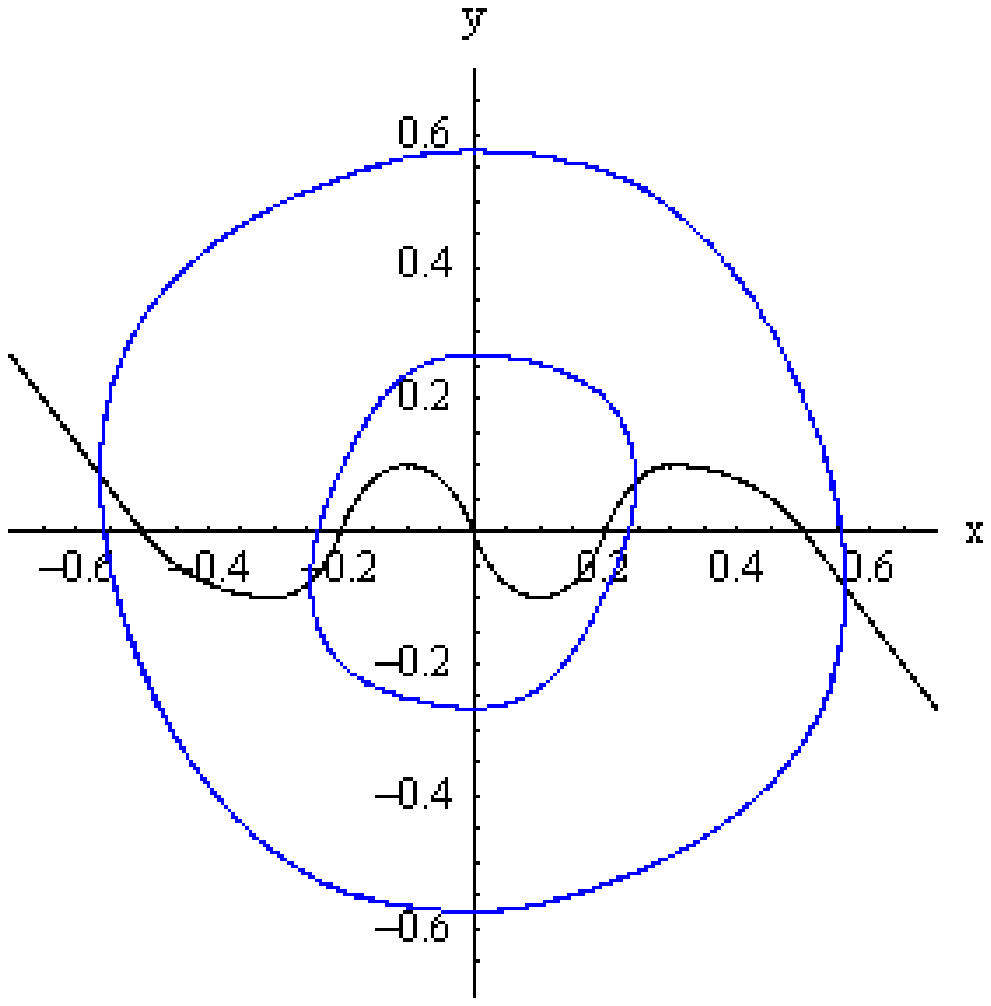}
\end{center}
\caption{Limit cycles for the system in Example \ref{Ex Construction 1} with
the curve of $F\left(  x\right)  .$}%
\label{Ex Construction 1 Fig}%
\end{figure}
\end{example}

\begin{remark}
\label{Choice Function Remark}From condition $\left(  C2\right)  $ in
$\cite{Odani N}$ we see that%
\[
g\left(  \phi_{k}\left(  x\right)  \right)  \phi_{k}^{\prime}\left(  x\right)
\geq g\left(  x\right)  \text{.}%
\]
If $g\left(  x\right)  =x$, then it gives
\[
\phi_{k}\left(  x\right)  \phi_{k}^{\prime}\left(  x\right)  \geq x\text{.}%
\]
Thus, in example $\ref{Ex Construction 1}$ if we take $\phi_{k}\left(
s\right)  =\phi_{1L}\left(  s\right)  =\sqrt{As^{2}+B}$, then the above
inequality gives%
\[%
\begin{array}
[c]{cc}
& \sqrt{As^{2}+B}\cdot\dfrac{2As}{2\sqrt{As^{2}+B}}\geq s\\
\text{i.e., } & As\geq s\\
\text{i.e., } & A\geq1\text{.}%
\end{array}
\]
By the definition of $f_{2L}^{\ast}$ and $f_{2}$ it follows that the remaining
part of the condition $\left(  C2\right)  $ is satisfied if%
\[
\left\vert F\left(  \phi_{1L}\left(  s\right)  \right)  \right\vert
\geq\left\vert F\left(  s\right)  \right\vert ,\quad s\in\left[  a_{0}%
,L_{0}\right]
\]
and in particular%
\[
\left\vert F\left(  \phi_{1L}\left(  0\right)  \right)  \right\vert
=\left\vert F\left(  0\right)  \right\vert \text{.}%
\]
Since, $\phi_{1L}\left(  s\right)  \in\left[  a_{1},a_{2}\right]  $ and
$s\in\left[  a_{0},L_{0}\right]  $ so it gives%
\begin{align}
\left\vert f_{2}\left(  \phi_{1L}\left(  s\right)  \right)  \right\vert  &
\geq\left\vert f_{1}\left(  s\right)  \right\vert ,\quad s\in\left[
a_{0},L_{0}\right] \nonumber\\
\text{i.e., }\left\vert H_{1L}\left(  f_{1L}\left(  s\right)  \right)
\right\vert  &  \geq\left\vert f_{1L}\left(  s\right)  \right\vert ,\quad
s\in\left[  a_{0},L_{0}\right]  \text{.} \label{Odani Condition}%
\end{align}
Next, in particular the equality occurs at $s=a_{0}=0$ and so we have%
\begin{align}
&  \left.  \hspace{0.32in}\left\vert F\left(  \phi_{1L}\left(  0\right)
\right)  \right\vert =\left\vert F\left(  0\right)  \right\vert \right.
\nonumber\\
&  \implies\left\vert H_{1L}\left(  f_{1L}\left(  0\right)  \right)
\right\vert =0\nonumber\\
&  \implies\left\vert H_{1L}\left(  0\right)  \right\vert =0\nonumber\\
&  \implies H_{1L}\left(  0\right)  =0 \label{Odani Condition Equality}%
\end{align}
since $F\left(  0\right)  =0$ and $f_{1L}\left(  0\right)  =0$. By our
construction we also see%
\[
s\cdot H_{1}\left(  s\right)  <0\quad\forall~s\text{.}%
\]
Thus, $\phi_{1L}$ behaves like choice function described by Odani. Here, the
condition $\left(  \ref{Odani Condition}\right)  $ does not hold for the
system discussed in example $\ref{Ex Construction 1}$. In fact, here%
\[
\left\vert H_{1L}\left(  f_{1L}\left(  s\right)  \right)  \right\vert
\leq\left\vert f_{1L}\left(  s\right)  \right\vert \quad s\in\left[
a_{0},L_{0}\right]  \text{.}%
\]
However, the conditions $($viz. $\bar{\alpha}_{i}<L_{i}$, etc.$)$ of theorem
$\ref{Th n Limit Cycle}$ are satisfied ensuring the existence of exactly two
limit cycles. This shows that the theorem $\ref{Th n Limit Cycle}$ and the
construction presented above covers a larger class of functions $F$ than those
covered in $\cite{Odani N}$. The equality in $\left(  \ref{Odani Condition}%
\right)  $ occurs in example $\ref{Ex Construction 1}$ only at the point
$s=a_{0}=0$. However, the equality can occur at points where $s\neq a_{0}$. We
present example $\ref{Ex Construction 2}$ below to show this kind of behaviour.
\end{remark}

\begin{remark}
The function $f_{2}$ in example $\ref{Ex Construction 1}$ is obtained from
$f_{1}$ by reflection and translation along $x$-axis. However, it is clear
from the construction of section $\ref{Construction}$, there is a plenty of
freedom in the possible extensions of $f_{1}$ having a fixed number of limit
cycles, as illustrated in examples $\ref{Ex Construction 2}$ and
$\ref{Ex Construction 3}$. In these examples we consider more general
transformations so that the limit cycles are obtained having amplitudes close
to those expected from the given physical $($dynamical$)$ problem.
\end{remark}

\begin{example}
\label{Ex Construction 2}In this example our target is to construct an example
in which%
\[
F_{+}\left(  x\right)  =\left\{
\begin{array}
[c]{l}%
0.055518-0.08\sqrt{1-\dfrac{\left(  x-0.144\right)  ^{2}}{0.04}},\quad0\leq
x\leq0.144\\
0.148506-0.172988\sqrt{1-\dfrac{\left(  x-0.144\right)  ^{2}}{\left(
0.206686\right)  ^{2}}},\quad0.144<x\leq0.34\\
0.0910146+0.0209854\sqrt{1-\dfrac{\left(  x-0.407\right)  ^{2}}{\left(
0.06751554\right)  ^{2}}},\quad0.34<x\leq0.407\\
-0.2280727+0.340073\sqrt{1-\dfrac{\left(  x-0.407\right)  ^{2}}{\left(
0.125376\right)  ^{2}}},\quad0.407<x\leq0.5\\
-3.0000372\left(  x-0.5\right)  ,\quad x>0.5\text{.}%
\end{array}
\right.
\]
and%
\[
F\left(  x\right)  =\left\{
\begin{array}
[c]{c}%
F_{+}\left(  x\right)  ,\quad x\geq0\\
F_{-}\left(  x\right)  ,\quad x<0
\end{array}
\right.
\]
Here, $a_{1}=0.2$, $a_{2}=0.5$, $L_{0}=0.144$ and $L_{1}=0.407$. It is easy to
show that $\phi_{1L}\left(  s\right)  =\sqrt{4.974392361\cdot s^{2}+0.0625}$
and $\phi_{1R}\left(  s\right)  =\sqrt{2.019706\cdot s^{2}+0.12376838}$. Here,%
\begin{align*}
f_{1}\left(  x\right)   &  =\left\{
\begin{array}
[c]{l}%
0.055518-0.08\sqrt{1-\dfrac{\left(  x-0.144\right)  ^{2}}{0.04}},\quad0\leq
x\leq0.144\\
0.148506-0.172988\sqrt{1-\dfrac{\left(  x-0.144\right)  ^{2}}{\left(
0.206686\right)  ^{2}}},\quad0.144<x\leq0.2
\end{array}
\right.  \\
f_{2}\left(  x\right)   &  =\left\{
\begin{array}
[c]{l}%
0.148506-0.172988\sqrt{1-\dfrac{\left(  x-0.144\right)  ^{2}}{\left(
0.206686\right)  ^{2}}},\quad0.2<x\leq0.34\\
0.0910146+0.0209854\sqrt{1-\dfrac{\left(  x-0.407\right)  ^{2}}{\left(
0.06751554\right)  ^{2}}},\quad0.34<x\leq0.407\\
-0.2280727+0.340073\sqrt{1-\dfrac{\left(  x-0.407\right)  ^{2}}{\left(
0.125376\right)  ^{2}}},\quad0.407<x\leq0.5
\end{array}
\right.
\end{align*}
The second part of the condition $\left(  C2\right)  $ in $\cite{Odani N}$
i.e., the condition $\left(  \ref{Odani Condition}\right)  $ does not hold. In
fact,
\begin{align*}
\left\vert H_{1L}\left(  f_{1L}\left(  s\right)  \right)  \right\vert  &
<\left\vert f_{1L}\left(  s\right)  \right\vert \text{ in }\left(
0,0.05290111\right)  \\
\text{and }\left\vert H_{1L}\left(  f_{1L}\left(  s\right)  \right)
\right\vert  &  >\left\vert f_{1L}\left(  s\right)  \right\vert \text{ in
}\left(  0.05290111,0.144\right)  \text{.}%
\end{align*}
The equality occurs at $s=0$ and $s=0.05290111$. Here, we get two limit cycles
crossing the $+ve$ $y$-axis at the points $\left(  0,0.29039755\right)  $ and
$\left(  0,0.567249\right)  $ respectively so that $y_{+}\left(  0\right)
=y_{-}\left(  0\right)  =0.29039755$ and $\bar{\alpha}_{1}=0.2892792083$.
Consequently, $\bar{\alpha}_{1}\leq L_{1}$ and the other conditions of theorem
$\ref{Th n Limit Cycle}$ with $g\left(  x\right)  =x$ are satisfied justifying
the existence of exactly two limit cycles. These two limit cycles alongwith
the curve of $F\left(  x\right)  $ has been shown in figure
$\ref{Ex Construction 2 Fig}$.\begin{figure}[h]
\begin{center}
\includegraphics[height=6cm]{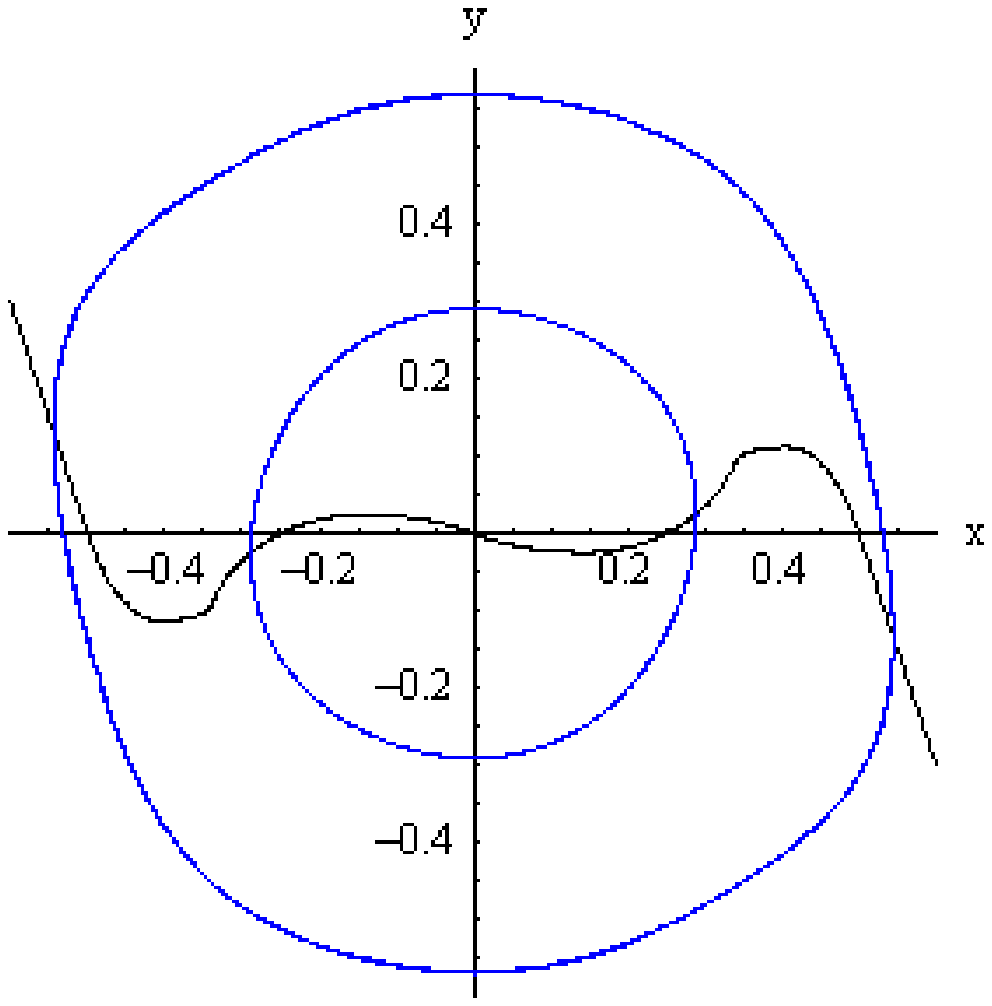}
\end{center}
\caption{Limit cycles for the system in Example \ref{Ex Construction 2} with
the curve of $F\left(  x\right)  .$}%
\label{Ex Construction 2 Fig}%
\end{figure}
\end{example}

\begin{example}
\label{Ex Construction 3}We now consider an example involving three limit
cycles by taking $a_{1}=0.1$, $a_{2}=0.2$, $a_{3}=0.4$ and%
\[
f_{1}\left(  x\right)  =0.04422166-0.08\sqrt{1-\frac{\left(  x-0.05\right)
^{2}}{\left(  0.06\right)  ^{2}}},\quad-0.1\leq x\leq0.1\text{.}%
\]
Here $L_{0}=0.05$ and let $L_{1}=0.15$, $L_{2}=0.3$. We take%
\begin{align*}
f_{2L}^{\ast}\left(  x\right)   &  =-0.04422166+0.08\sqrt{1-\frac{\left(
x-0.15\right)  ^{2}}{\left(  0.06\right)  ^{2}}}\\
f_{2R}^{\ast}\left(  x\right)   &  =-0.04422166+0.08\sqrt{1-\frac{\left(
x-0.15\right)  ^{2}}{\left(  0.06\right)  ^{2}}}\text{.}%
\end{align*}
It is easy to construct%
\begin{align*}
\phi_{1L}\left(  s\right)   &  =\sqrt{5s^{2}+0.01},\quad s\in\left[
a_{0},L_{0}\right]  \\
\phi_{1R}\left(  s\right)   &  =\sqrt{\frac{7}{3}s^{2}+\frac{5}{300}},\quad
s\in\left[  L_{0},a_{1}\right]  \text{.}%
\end{align*}
Next, we take%
\begin{align*}
f_{3L}^{\ast}\left(  x\right)   &  =0.0043819183-0.03\sqrt{1-\frac{\left(
x-0.3\right)  ^{2}}{\left(  0.101084111\right)  ^{2}}}\\
f_{3R}^{\ast}\left(  x\right)   &  =0.0043819183-0.03\sqrt{1-\frac{\left(
x-0.3\right)  ^{2}}{\left(  0.101084111\right)  ^{2}}}\text{.}%
\end{align*}
We can similarly construct%
\begin{align*}
\phi_{2L}\left(  s\right)   &  =\sqrt{4s^{2}}=2s,\quad s\in\left[  a_{1}%
,L_{1}\right]  \\
\phi_{2R}\left(  s\right)   &  =2s,\quad s\in\left[  L_{1},a_{2}\right]
\text{.}%
\end{align*}
so that $\phi_{2L}\left(  a_{1}\right)  =a_{2}$, $\phi_{2L}\left(
L_{1}\right)  =L_{2}$, $\phi_{2R}\left(  L_{1}\right)  =L_{2}$ and $\phi
_{2R}\left(  a_{2}\right)  =a_{3}$. We define%
\[
F_{+}\left(  x\right)  =\left\{
\begin{array}
[c]{l}%
0.04422166-0.08\sqrt{1-\dfrac{\left(  x-0.05\right)  ^{2}}{\left(
0.06\right)  ^{2}}},\quad0\leq x<0.1\\
-0.04422166+0.08\sqrt{1-\dfrac{\left(  x-0.15\right)  ^{2}}{\left(
0.06\right)  ^{2}}},\quad0.1\leq x<0.2\\
0.0043819183-0.03\sqrt{1-\dfrac{\left(  x-0.3\right)  ^{2}}{\left(
0.101084111\right)  ^{2}}},\quad0.2\leq x<0.4\\
2.0100758\left(  x-0.4\right)  ,\quad0.4\leq x\text{.}%
\end{array}
\right.
\]%
\[
F\left(  x\right)  =\left\{
\begin{array}
[c]{c}%
F_{+}\left(  x\right)  ,\quad x\geq0\\
F_{-}\left(  x\right)  ,\quad x<0
\end{array}
\right.
\]
to make $F$ continuously differentiable. We can easily calculate that
$\bar{\alpha}_{1}=0.12418214965$ and $\bar{\alpha}_{2}=0.2354818163$ and
consequently $\bar{\alpha}_{i}<L_{i}$ for $i=1,2$. All the other conditions of
theorem $\ref{Th n Limit Cycle}$ with $g\left(  x\right)  =x$ are satisfied
and hence we get three distinct limit cycles as shown in figure
$\ref{Ex Construction 3 Fig}$ alongwith the curve of $F\left(  x\right)  $
defined above.\begin{figure}[h]
\begin{center}
\includegraphics[height=6cm]{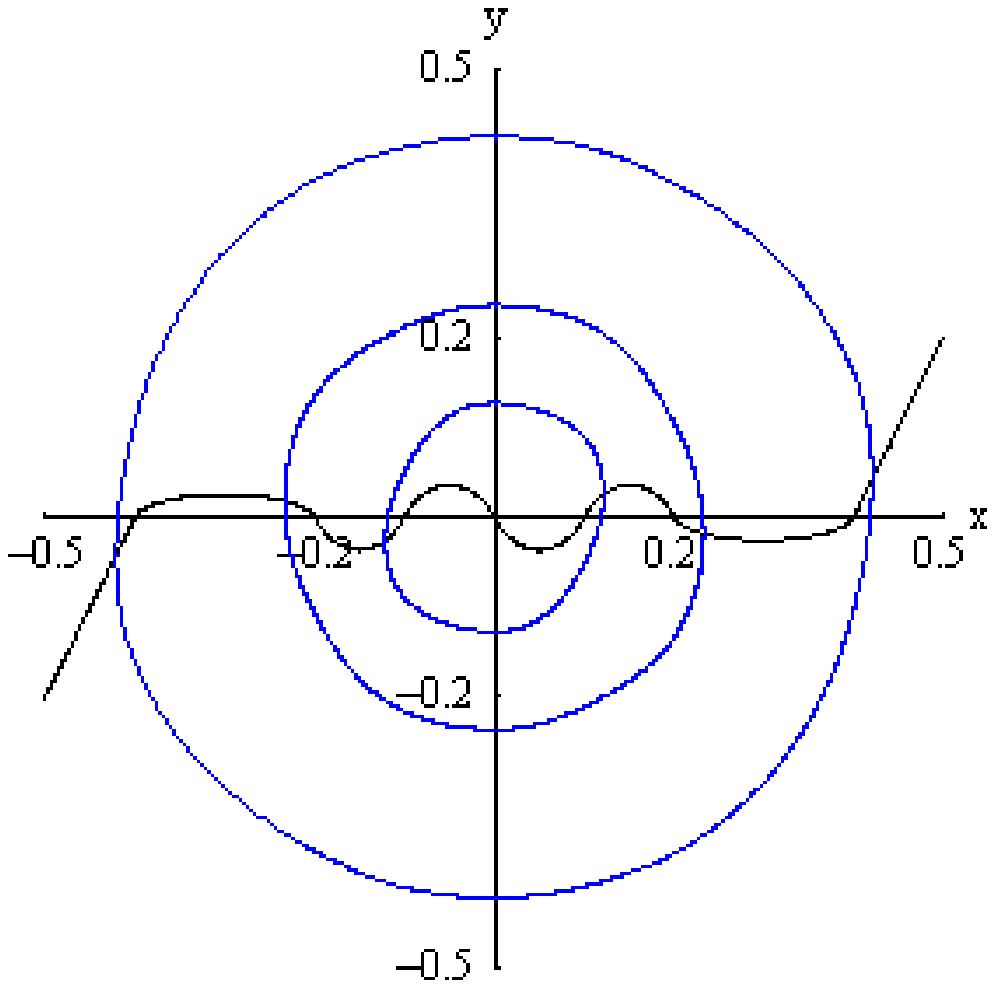}
\end{center}
\caption{Three limit cycles for the system in Example \ref{Ex Construction 3}
with the curve of $F\left(  x\right)  .$}%
\label{Ex Construction 3 Fig}%
\end{figure}
\end{example}

\begin{remark}
Here the function $F$ is defined in such a manner that $\left\vert F\left(
L_{0}\right)  \right\vert >\left\vert F\left(  L_{2}\right)  \right\vert $
implying that $\beta_{2}$ mentioned in Theorem 3 of $\cite{Chen Llibre Zhang}$
or in theorem 7.12, chapter 4 of the book $\cite{Zhing Tongren Wenzao}$ does
not exist and hence these theorems are not applicable for the corresponding
Lienard system.
\end{remark}

\end{document}